\documentclass[12pt]{article}  
\usepackage{amsmath}
\usepackage{theorem}                 
\usepackage{amssymb}
\usepackage{xypic}

\newtheorem{theorem}{Theorem}
\newtheorem{corollary}{Corollary}
\newtheorem{proposition}{Proposition}

\theorembodyfont{\rmfamily}
\newtheorem{example}{Example}

\newenvironment{examples}
{\smallskip\noindent{\bf Examples\/}.}{\smallskip\par}
\newenvironment{remark}
{\smallskip\noindent{\bf Remark\/}.}{\smallskip\par}
\newenvironment{remarks}
{\smallskip\noindent{\bf Remarks\/}.}{\smallskip\par}
\newenvironment{proof}
{\noindent{\it Proof\/}.}{{ \hfill $\Box$}\smallskip\par}

\newcommand{\CC}{{\Bbb C}}

\newcommand{\RR}{{\Bbb R}}

\newcommand{\calO}{{\cal O}}

\newcommand{\calA}{{\cal A}}

\newcommand{\calJ}{{\cal J}}

\newcommand{\eps}{\varepsilon}
\newcommand{\icis}{{ICIS}}
\newcommand{\ELKh}{{\scriptscriptstyle \rm ELKh}}

\title{Quadratic forms for a 1-form on an isolated complete
intersection singularity}
\author{W.~Ebeling and S.~M.~Gusein-Zade
\thanks{Partially supported by the DFG-programme ''Global methods in
complex geometry'' (Eb 102/4--3), grants RFBR--04--01--00762,
NSh--1972.2003.1.
}
}

\date{}

\begin{document}

\maketitle

\begin{abstract}
We consider a holomorphic 1-form $\omega$ with an isolated zero on an isolated complete intersection
singularity $(V,0)$. We construct quadratic forms on an algebra of functions and on a module of
differential forms associated to the pair $(V,\omega)$. They generalize the 
Eisenbud--Levine--Khimshiashvili quadratic form defined for a smooth $V$.
\end{abstract}

\section*{Introduction}

Let $g=(g_1, \ldots, g_n)$ be the germ of an analytic map $(\CC^n,
0)\to(\CC^n, 0)$
which is real (i.e. $g(\RR^n)\subset\RR^n$) and finite. The latter fact
is equivalent to
the one that the local algebra $\calA_g=\calO_{\CC^n, 0}/\langle g_1, \ldots, g_n \rangle$ is
a finite dimensional vector space (here $\calO_{\CC^n, 0}$ is the ring of germs
at the origin of analytic functions on $\CC^n$, $\langle g_1, \ldots, g_n \rangle$ is the
ideal generated by the corresponding germs). The degree $\nu$ of the map
$g:(\CC^n, 0)\to(\CC^n, 0)$ is equal to the dimension $\dim_\CC \calA_g$ of the
algebra $\calA_g$ (as a $\CC$-vector space). Let $\calA_g^\RR$ be the real part
of the algebra $\calA_g$ (it is the ring of germs at $0$ of real analytic
functions on $\RR^n$ factorized by the ideal
$\langle g_1, \ldots, g_n \rangle$; $\dim_\RR \calA_g^\RR=\dim_\CC \calA_g$). The famous
theorem
of D.~Eisenbud, H.~Levine and G.~Khimshiashvili states that there exists a
non-degenerate
quadratic form $Q_{\ELKh}$ on the algebra $\calA_g$ such that it is real on
$\calA_g^\RR$
and the degree of the real analytic map
$g_\RR=(g_1, \ldots, g_n):(\RR^n, 0)\to(\RR^n, 0)$ is equal to the signature of
the quadratic form $Q_{\ELKh}$ on the algebra $\calA_g^\RR$. Moreover this
quadratic form
can be defined by the formula
$$
Q_{\ELKh}(\varphi,\psi)=R(\varphi\cdot \psi)
$$
where $R$ is an arbitrary real linear function on the space $\calA_g$
which has
a positive value on the Jacobian
$$
{\rm Jac}_g={\rm Jac}(g_1, \ldots, g_n)=\det\left( \frac{\partial g_i}{\partial x_j} \right)
$$
of the map $g$.

A linear function with this property can also be defined by the formula
$$
R(\varphi)=\lim\limits_{\eps\to 0}\ \sum\limits_{i=1}^\nu
\frac{\varphi(P_i)}{{\rm Jac}_g(P_i)}
$$
where the sum is over all points $P_i$ of the preimage
$g^{-1}(\eps)\subset\CC^n$
of a regular value $\eps\in\CC^n$ of the map $g$ (cf. \cite[\S 5]{AGV}). The
Eisenbud--Levine--Khimshiashvili theorem can be also formulated in terms of
indices of vector fields or of 1-forms on $(\RR^n,0)$.

There were some results which can be considered as attempts to generalize the
Eisenbud--Levine--Khimshiashvili theorem to singular varieties: see, e.g.,
\cite{MvS, GMM1, GMM2, Sz, EG1}.

Let $f=(f_1, \ldots, f_k):(\CC^n, 0)\to(\CC^k, 0)$ be a germ of a complex
analytic map
which  defines an isolated complete
intersection singularity (\icis) $(V,0)=f^{-1}(0)\subset(\CC^n,0)$ and let 
$\omega$ be
the germ of a holomorphic 1-form on $(\CC^n,0)$ whose restriction to $V$ has
no singular
points (zeroes) in a punctured neighbourhood of the origin. 
Let $I_{V;\omega}$ be the ideal of the ring $\calO_{\CC^n,0}$ generated by
the functions $f_1$, \dots, $f_k$ and by the $(k+1)\times(k+1)$-minors of
the matrix
$$
\left(
\begin{array}{ccc}
\frac{\partial f_1}{\partial x_1} & \cdots &
\frac{\partial f_1}{\partial x_n} \\
 \vdots & \ddots & \vdots \\
\frac{\partial f_k}{\partial x_1} & \cdots &
\frac{\partial f_k}{\partial x_n} \\
A_1 & \cdots & A_n
\end{array}
\right) 
$$
and let
$\calA_{V;\omega}$ be the factor algebra $\calO_{\CC^n,0}/I_{V;\omega}$.
For $\omega = dh$, $h \in {\cal O}_{\CC^n,0}$, this algebra was considered in \cite {Gr1, Le} in
connection with the computation of the Milnor number of an \icis.

Suppose that the map $f$ and the 1-form $\omega$ are real 
(i.e. $f(\RR^n)\subset\RR^k$ and
$\omega$ has real values on the tangent spaces $T_x\RR^n$).
The 1-form
$\omega$ considered
on the real part $(V_\RR, 0)$ of the {\icis} $(V,0)$ has an invariant~--- an
index
${\rm ind}_{V,0}\,\omega$ (which can be called the {\em radial index}): \cite{EG1}. 
In
\cite{EG1},
there was constructed a family $Q_\eps$ of quadratic forms on a vector space (the dimension of which
coincides with the dimension of the algebra $\calA_{V;\omega}$ as a $\CC$-vector space) 
depending on $\eps$ from a neighbourhood of the origin in the
target
$\CC^k$ of the map $f$ such that:\newline
(1) the quadratic form $Q_\eps$ is nondegenerate for $\eps$ from the complement
to the bifurcation diagram $\Sigma\subset(\CC^k,0)$ of the map $f$;\newline
(2) for real $\eps$, the quadratic form
$Q_{\eps}$ is real and, for real $\eps\notin\Sigma$, its signature is equal to
$$
{\rm ind}_{V,0}\, \omega + (\chi(V_{\eps})-1),
$$
where $V_\eps=f^{-1}(\eps) \cap \RR^n \cap B_\delta$ for $0<\eps \ll \delta$ small
enough,
$B_\delta$ is the ball of radius $\delta$ centred at the origin in $\CC^n$.

Here we show that the construction of \cite{EG1} gives a (non-trivial)
quadratic
form on the algebra $\calA_{V;\omega}$. For the curve case (i.e. for
$k=n-1$) this
quadratic form is nondegenerate and coincides with the one constructed by
J.~Montaldi
and D.~van~Straten in \cite{MvS}.

We also define a quadratic form on a quotient module $\Omega_{V;\omega}$
of the module $\Omega_{\CC^n,0}^{n-k}$ of $(n-k)$-forms on $(\CC^n,0)$ which,
as a $\CC$-vector space, has the same dimension as $\calA_{V;\omega}$.
We describe some relations between these quadratic forms.

We define these quadratic forms by an analytic construction. It would be
interesting
to understand to which extend it is possible to define them in algebraic terms. 

\section{Jacobian of a 1-form}

Let $M^m$ be a complex analytic manifold of dimension $m$ with local
coordinates
$y^1$, \dots, $y^m$. A vector field
$X=\sum\limits_{i=1}^m X^i\frac{\partial}{\partial y^i}$
on $M$ is a tensor field of type $(1, 0)$. Its "derivative" (Jacobian matrix)
$\left(\frac{\partial X^i}{\partial y^j}\right)$ is generally speaking not
a tensor,
but it is a tensor of type $(1, 1)$ at points where the vector field $X$
vanishes.
At such a point $P$ it has a well defined (i.e. not dependent on the choice
of coordinates) determinant which is the value of the Jacobian
${\rm Jac}(X^1, \ldots, X^m)$ of the vector field $X$ at the point $P$.

Now let $\omega=\sum\limits_{i=1}^m A_i\,dy^i$ be a 1-form on $M$,
i.e. a tensor field of type $(0, 1)$. Again,
its Jacobian matrix $\calJ=\left(\frac{\partial A_i}{\partial y^j}\right)$
defines the tensor
$\sum\limits_{i,j}\frac{\partial A_i}{\partial y^j}\, dy^i\otimes dy^j$
of type $(0, 2)$ (i.e. a bilinear form on the tangent space) only at points
where the 1-form $\omega$ vanishes. The determinant
$J$ of the matrix $\calJ=\left(\frac{\partial A_i}{\partial y^j}\right)$ is not
a scalar (it depends on the choice of local coordinates). Under a change of
coordinates it is multiplied by the square of the Jacobian of the change
of coordinates (since the Jacobian matrix $\calJ$ is transformed to
$C^T\calJ C$).
Therefore it should be considered as the coefficient in the tensor
$$
J(dy^1\wedge\ldots\wedge dy^m)^{\otimes 2}
$$
of type $(0,2m)$. In this sense the Jacobian of a 1-form is a sort of a
"quadratic differential". To get a number, one can divide this tensor by the
tensor square of a volume form.

Let $m=n-k$ and let $M$ be a submanifold in $\CC^n$ defined by $k$ equations,
i.e. let $M=F^{-1}(0)$, where $F=(f_1, \ldots, f_k)$ is a nondegenerate
holomorphic map from $\CC^n$ to $\CC^k$ (defined in a neighbourhood of the
manifold $M$). Let us fix volume forms on $\CC^n$ and $\CC^k$, say, the
standard ones
$\sigma_n=dx_1\wedge\ldots\wedge dx_n$ and $\sigma_k=dz_1\wedge\ldots\wedge
dz_k$
where $x_1$, \dots, $x_n$ and $z_1$, \dots, $z_k$ are Cartesian coordinates
in $\CC^n$
and in $\CC^k$ respectively. There exists (at least locally) an
$(n-k)$-form $\sigma$
on $\CC^n$ such that $\sigma_n=F^*\sigma_k\wedge\sigma$. The restriction
$\sigma_{\vert M}$ of the form $\sigma$ to the manifold $M$ is well defined
and is
a volume form on $M$. We also denote it by $\sigma$. After that, for a
holomorphic
1-form $\omega$
on the manifold $M$ (say, the restriction to $M$ of a 1-form on the ambient space
$\CC^n$),
the above construction gives well defined numbers $\widetilde J(P)$
associated with
singular points (zeroes) $P$ of the form $\omega$ ("the values of the
Jacobian").

Let $\omega=\sum\limits_{i=1}^n A_i\, dx_i$ be a holomorphic 1-form on
$\CC^n$ and
let $P$ be a singular point (zero) of the form $\omega$ on $M$. There exists
a $(k\times k)$-minor of the matrix $\left(\frac{\partial f_i}{\partial
x_j}\right)$
which is different from zero at the point $P$. Let it be
$$
\Delta=\left\vert\frac{\partial f_i}{\partial x_j}\right\vert_{i,j=1,
\ldots, k}.
$$
In this case $x_{k+1}$, \dots, $x_n$ are local coordinates on $M$ in a
neighbourhood
of the point $P$ and therefore the restriction of the 1-form $\omega$ to
the manifold
$M$ can be written as
$$
\sum\limits_{i=k+1}^n \widehat A_i\, dx_i.
$$
The Jacobian determinant of the 1-form $\omega$ in these coordinates is
$$
J=\left\vert\frac{\partial\widehat A_i}{\partial x_j}\right\vert_{i,j=k+1,
\ldots, n}.
$$
A precise formula for the determinant $J$ can be found in \cite[Prop.~4]{EG1}.
For $i=k+1, \dots, n$, let
$$
m_{i}
:= \left|
\begin{array}{cccc}
\frac{\partial f_1}{\partial x_1} &
\cdots &
\frac{\partial f_1}{\partial x_k} &
\frac{\partial f_1}{\partial x_i} \\
\vdots & \ddots & \vdots & \vdots \\
\frac{\partial f_k}{\partial x_1} &
\cdots &
\frac{\partial f_k}{\partial x_k} &
\frac{\partial f_k}{\partial x_i}
\\ A_1 & \cdots & A_k & A_i
\end{array}
\right|.
$$
One can see that
$$
\omega_{\vert M} =  \frac{m_{k+1}}{\Delta} dx_{k+1} +
\cdots + \frac{m_n}{\Delta}dx_n.
$$
{}From this it was derived that
$$
J =
\frac{1}{\Delta^{2+(n-k)}}
\left|
\begin{array}{cccc}
\Delta & \frac{\partial \Delta}{\partial x_1} & \cdots &
\frac{\partial\Delta}{\partial x_n}\\
0 & \frac{\partial f_1}{\partial x_1} & \cdots &
\frac{\partial f_1}{\partial x_n} \\
\vdots & \vdots & \ddots & \vdots \\
0 & \frac{\partial f_k}{\partial x_1} & \cdots &
\frac{\partial f_k}{\partial x_n} \\
m_{k+1} & \frac{\partial m_{k+1}}{\partial x_1} & \cdots &
\frac{\partial m_{k+1}}{\partial x_n}\\
\vdots & \vdots & \ddots & \vdots \\
m_n & \frac{\partial m_n}{\partial x_1} & \cdots &
\frac{\partial m_n}{\partial x_n}
\end{array}
\right|.
$$
At the point $P$ (a zero of the 1-form $\omega$ on the manifold $M$) the
minors $m_i$
vanish and therefore
\begin{eqnarray*}
J(P) & = &
\frac{1}{\Delta^{1+(n-k)}}
\left|
\begin{array}{ccc}
\frac{\partial f_1}{\partial x_1} & \cdots &
\frac{\partial f_1}{\partial x_n} \\
 \vdots & \ddots & \vdots \\
\frac{\partial f_k}{\partial x_1} & \cdots &
\frac{\partial f_k}{\partial x_n} \\
\frac{\partial m_{k+1}}{\partial x_1} & \cdots &
\frac{\partial m_{k+1}}{\partial x_n}\\
\vdots & \ddots & \vdots \\
\frac{\partial m_n}{\partial x_1} & \cdots &
\frac{\partial m_n}{\partial x_n}
\end{array}
\right| \\
& = & \frac{1}{\Delta^{1+(n-k)}}{\rm Jac}(f_1, \ldots, f_k, m_{k+1}, \ldots,
m_n)
\end{eqnarray*}
(the expressions have to be evaluated at the point $P$ of course).

One has
$$
F^*\sigma_k=\sum\limits_{1\le s_1<\ldots<s_k\le n}
\left\vert\frac{\partial f_i}{\partial x_{s_j}}\right\vert_{i,j=1, \ldots, k}
dx_{s_1}\wedge\ldots\wedge dx_{s_k}.
$$
Therefore $\sigma=\frac{1}{\Delta}dx_{k+1}\wedge\ldots\wedge dx_n$ and
$$
\widetilde{J}(P)=\Delta^2 J(P).
$$

\begin{remark}
The number $\widetilde J(P)$ was already used in \cite{EG1}.
\end{remark}

\section{A quadratic form on the algebra of functions.}

Let $f=(f_1, \ldots, f_k):(\CC^n,0)\to(\CC^k,0)$ be a holomorphic map which
defines an
$(n-k)$-dimensional isolated complete intersection singularity (\icis)
$V=f^{-1}(0)\subset(\CC^n,0)$ and let $\omega=\sum\limits_{i=1}^n A_i dx_i$
be the germ of a 1-form on $(\CC^n, 0)$. Let $F:(\CC^n\times\CC_\eps^M,0)\to
(\CC^k\times\CC_\eps^M,0)$ be a deformation of the map
$f:(\CC^n,0)\to(\CC^k,0)$
($F(x, \eps)=(f_\eps(x), \eps)$, $f_0=f$, $f_\eps=(f_{1\eps}, \ldots,
f_{k\eps})$) and let $\omega_\eps$ be a deformation of the form $\omega$
(defined in a neighbourhood of the origin in $\CC^n\times\CC_\eps^M$) such
that, for generic
$\eps\in (\CC_\eps^M,0)$, the preimage $f_\eps^{-1}(0)$ is smooth and
(the restriction of) the form $\omega_\eps$ to it has only non degenerate
singular
points.

The number $\nu$ of these points is the (complex) index of the 1-form $\omega$
on $(V,0)$ considered in \cite{EG1} (it is an analogue of the GSV-index of a
vector field).
It is equal to the dimension (as a $\CC$-vector space) of the algebra
$\calA_{V;\omega}$
(\cite{EG1}; see also \cite{EG2} where some inaccuracies in the proof were
corrected).

Let $\Sigma\subset (\CC^M_\eps,0)$ be the germ of the set of the values of
the parameters $\eps$ from $(\CC^M_\eps,0)$ such that either the preimage
$f_\eps^{-1}(0)$ is singular or the restriction of the 1-form $\omega_\eps$
to it has degenerate singular points. The {\em bifurcation diagram} $\Sigma$
is a germ of a hypersurface in $(\CC^M_\eps, 0)$. For $\eps\notin\Sigma$,
let $P_1$, \dots, $P_\nu$ be the (non-degenerate) singular points of the 1-form
$\omega_\eps$ on the $(n-k)$-dimensional manifold $f_\eps^{-1}(0)$.
For $\eps\notin\Sigma$ and a germ $\varphi\in\calO_{\CC^n,0}$, let
\begin{equation}\label{Formula1}
R(\varphi,\eps):=\sum\limits_{i=1}^\nu \frac{\varphi(P_i)}{\widetilde J_\eps(P_i)}.
\end{equation}
For a fixed $\varphi$ the function $R_\varphi(\eps):=R(\varphi,\eps)$ is holomorphic in
the complement of the bifurcation diagram $\Sigma$.

\begin{theorem} \label{Thm1}
The function $R_\varphi(\eps)$ has removable singularities on the bifurcation
diagram $\Sigma$.
\end{theorem}

\begin{proof} 
Without any loss of generality we can suppose that the deformation $(F,\omega_\eps)$ is ''as big as we
would like''. In particular we can assume that it includes a versal deformation of the map $f$ with
the trivial deformation of the 1-form $\omega$ and also that $n$ parameters of it (marked by elements
of $\CC^{n \ast}$) correspond to the trivial deformation of the map and adding the corresponding linear
form (i.e.\ a 1-form with constant coefficients) to the 1-form.

It is sufficient to prove that the function $R_\varphi(\eps)$ has removable
singularities outside of a set of codimension 2 in $\CC^M_\eps$. Therefore it
suffices to prove this for points of the bifurcation diagram $\Sigma$ of two
types:
\begin{enumerate}
\item  Those points $\eps \in \Sigma$ for which the preimage $f^{-1}_\eps(0)$ is
smooth (but a singular point of the 1-form $\omega_\eps$ on it is degenerate).
In this case the proof can be found, e.g., in \cite[\S 5]{AGV}.
\item Those points $\eps \in \Sigma$ for which the preimage $f^{-1}_\eps(0)$ has
one singular point $P$ of type $A_1$, the 1-form $\omega_\eps$ does not vanish
at this point in $\CC^n$, and its kernel is not tangent to the tangent cone of
the variety $f^{-1}_\eps(0)$. In this case the proof is essentially already
contained in \cite{EG1}. For the sake of completeness we repeat it in an
appropriate way.
\end{enumerate}

Without loss of generality (excluding equations which are non-deg\-en\-er\-ate), we
can suppose that $k=1$, $n \geq 2$. Changing local coordinates in a
neighbourhood of the point $P$, we can suppose that $P$ is the origin in
$\CC^n$, $f_1=x_1^2+ \ldots + x_n^2$, $\omega(0)=dx_1$ (we omit $\eps$ in the
notations). The last equation means that $\omega=(1+C_1)dx_1+C_2dx_2 + \ldots
C_ndx_n$ where $C_i(0)=0$. It is sufficient to show that the function $R_\varphi$
has a finite limit when $\eta \to 0$ for the deformation
$f_{1\eta}=f_1-\eta^2$, $\eta \in \CC$, of the equation keeping the 1-form
$\omega$ constant. For $\eta \neq 0$ small enough the 1-form $\omega$ has two
singular points $P_+$ and $P_-$ on the level manifold $F_{1\eta}=0$ with
coordinates $x_1= \pm \eta + o(\eta)$, $x_i=o(\eta)$ for $i \geq 2$. The
formula (\ref{Formula1}) gives
$$\widetilde{J}(P_\pm)=(-1)^{n-1} 4 (\pm \eta)^{3-n} (1+o(\eta)).$$
The corresponding summands in the expression for $R_\varphi$ are
$$\frac{(-1)^{n-1}}{4} \left(
\varphi(P_+)(\eta^{n-3}+o(\eta^{n-3}))+\varphi(P_-)((-\eta)^{n-3}+o(\eta^{n-3}))
\right).$$
For $n \geq 3$ or for $\varphi(0)=0$, each summand has a finite limit when $\eta
\to 0$. For $n=2$ and $\varphi(x) \equiv 1$,
$\varphi(P_+)\eta^{-1}+\varphi(P_-)(-\eta)^{-1}=0$.
\end{proof}

Theorem~\ref{Thm1} means that $R_\varphi(\eps)$ can be extended to a holomorphic
function on $(\CC^M_\eps,0)$ (which we denote by the same symbol). Let
\begin{equation}
R(\varphi):=R_\varphi(0).
\end{equation}
This defines a linear function $R$ on ${\cal O}_{\CC^n,0}$. It is clear
that, if
the map $f$ and the 1-form $\omega$ are real, the function $R$ is real as well.

\begin{remark}
Just in the same way we can assume that the function $\varphi$ depends on the parameters $\eps \in
\CC^M_\eps$ as well: $\varphi= \varphi(x,\eps)$, $\varphi : (\CC^n \times \CC^M_\eps,0) \to (\CC,0)$.
The statement and the proof are the same. Moreover, the limit $R(\varphi)= \lim\limits_{\eps \to 0}
R_\varphi(\eps)$ is also the same, i.e.\ $R(\varphi(x,\eps))=R(\varphi(x,0))$ (this was the reason why
we did not include the statement into the theorem). To show the latter fact we can enlarge the
deformation space including one with the function $\varphi$ depending on $\eps$ and the identical one
with the function $\varphi$ not changing with $\eps$.
\end{remark}

As above let $I_{V;\omega} \subset {\cal O}_{\CC^n,0}$ be the ideal generated
by the equations $f_1$, \dots , $f_k$ of the {\icis} $(V,0)$ and by the $(k+1)
\times (k+1)$-minors of the matrix
\begin{equation}\begin{pmatrix} \frac{\partial f_1}{\partial x_1} & {\cdots } &
 \frac{\partial f_1}{\partial x_{n}} \\ {\vdots} & {\cdots} &
{\vdots} \\ \frac{\partial
f_k}{\partial x_1} & \cdots &
 \frac{\partial f_k}{\partial x_{n}} \\
A_1 & \cdots & A_{n}
\end{pmatrix}\,. \label{Matrix1}
\end{equation}

\begin{proposition} \label{Prop1}
The linear function $R$ vanishes on the ideal $I_{V;\omega}$.
\end{proposition}

\begin{proof} To prove that $R$ vanishes on a germ $\varphi=hf_i$ ($h \in {\cal
O}_{\CC^n,0}$), consider the function $\varphi(x,\eps)=h(f_i + \eps_i)$. On a
nonsingular level set $f+\eps=0$ the function $\varphi(x,\eps)$ vanishes
identically and therefore the expression (\ref{Formula1}) tends to zero. To
prove that $R$ vanishes on a germ $\varphi=hm_I$ ($h \in {\cal
O}_{\CC^n,0}$) where $m_I$ is one of the $(k+1) \times (k+1)$-minors of the
matrix (\ref{Matrix1}), consider the function $\varphi(x,\alpha)= h \cdot
m_I(x,\alpha)$ where $m_I(x,\alpha)$ is the corresponding minor of the matrix
\begin{equation}\begin{pmatrix} \frac{\partial f_1}{\partial x_1} & {\cdots } &
 \frac{\partial f_1}{\partial x_{n}} \\ {\vdots} & {\cdots} &
{\vdots} \\ \frac{\partial
f_k}{\partial x_1} & \cdots &
 \frac{\partial f_k}{\partial x_{n}} \\
A_1-\alpha_1 & \cdots & A_{n}-\alpha_n
\end{pmatrix}\,. \label{Matrix2}
\end{equation}
The function $\varphi(x,\alpha)$ vanishes at all the singular points of the 1-form
$\omega_\alpha$ on the level set $\{ f_\eps=0\}$ and therefore the expression
(\ref{Formula1}) tends to zero.
\end{proof}

\begin{corollary} The formula (\ref{Formula1}) defines a linear function $R$ on
the (finite dimensional) algebra ${\cal A}_{V;\omega} = {\cal
O}_{\CC^n,0}/I_{V;\omega}$.
\end{corollary}

\begin{remark} The function $R$ is well-defined up to multiplication by a unit
in the algebra ${\cal A}_{V;\omega}$ (which depends on the choice of volume
forms on $(\CC^n,0)$ and on $(\CC^k,0)$).
\end{remark}

At the moment, the definition of the linear function $R$ on ${\cal
A}_{V;\omega}$ uses the map $f$ (i.e. equations of the {\icis} $(V,0)$) and a
1-form $\omega$ defined on the ambient space $(\CC^n,0)$. It is clear that up
to multiplication by a unit the function $R$ depends only on the {\icis} $(V,0)$
itself. It fact it merely depends on the restriction of the 1-form $\omega$ to
the {\icis} $(V,0)$.

\begin{proposition} \label{Prop2}
The linear function $R$ is defined by the class of the 1-form $\omega$ in the
module
$$\Omega^1_{V,0} = \Omega^1_{\CC^n,0}/(f_i \Omega^1_{\CC^n,0}, {\cal
O}_{\CC^n,0} df_i)$$
of germs of 1-forms on the {\icis} $(V,0)$.
\end{proposition}

\begin{proof} Let
$$\omega'=\omega + \sum_{i=1}^k f_i \eta_i + \sum_{i=1}^k h_i df_i \quad
(\eta_i \in \Omega^1_{\CC^n,0}, \ h_i \in {\cal O}_{\CC^n,0}).$$
Consider the deformation
$$\omega'_\eps = \omega_\eps + \sum_{i=1}^k f_{i\eps} \eta_i + \sum_{i=1}^k h_i
df_{i\eps}$$
of the 1-form $\omega'$. On a (non-singular) level set $\{ f_{i\eps} =0\}$ the
1-form $\omega_\eps'$ is equal to $\omega_\eps$ and therefore the expressions (\ref{Formula1}) for
$\omega_\eps$ and $\omega'_\eps$ coincide.
\end{proof}

The linear function $R$ on the algebra ${\cal A}_{V;\omega}$ defines in the
natural way a quadratic form $Q^{\cal A}_{V;\omega}$ on it:
$$
Q^{\cal A}_{V;\omega}(\varphi, \psi):= R(\varphi \cdot \psi).
$$

\begin{example} \label{Ex1}
Let $k=1$, $f_1(x)=f(x)=\sum_{i=1}^n x_i^2$, $\omega = \sum_{i=1}^n a_i x_i
dx_i (= \frac{1}{2} d(\sum a_i x_i^2))$, where $a_i$ are pairwise different.
The algebra
$${\cal A}_{V;\omega} = {\cal O}_{\CC^n,0}/\left\langle \sum x_i^2, x_ix_j \right\rangle$$
has dimension $2n$ as a vector space over $\CC$. It is generated by the classes
of the monomials $1$, $x_1$, \dots , $x_n$, $x_1^2$, \dots , $x_{n-1}^2$. On
the level set $\{f=\eps^2\}$ the singular points of the 1-form $\omega$ are $P_i^\pm = (0, \ldots ,
0, \pm \eps , 0, \ldots ,0)$. They are non-degenerate.
The value of $\widetilde{J}$ at the point
$P_i^\pm$ is equal to $\eps^2 \prod_{j \neq i} (a_j-a_i)$. One has
$$R(1) = \lim_{\eps \to 0} \frac{2}{\eps^2} \left( \sum_{i=1}^n
\frac{1}{\prod_{j \neq i} (a_j-a_i)} \right) = 0$$
since the expression in parentheses is identically equal to zero; $R(x_i)=0$
(since the values of $x_i$ at the points $P^+_i$ and $P^-_i$ differ only by
their signs), and finally
$$R(x_i^2) = \frac{2}{\prod_{j \neq i} (a_j-a_i)}.$$
Therefore
$$Q^{\cal A}_{V;\omega}(1,x_i^2)=Q^{\cal
A}_{V;\omega}(x_i,x_i)=\frac{2}{\prod_{j \neq i} (a_j-a_i)}$$
and $Q^{\cal A}_{V;\omega}(\varphi, \psi)=0$ for all other pairs of the generators
listed above. Hence ${\rm rk}\, Q^{\cal A}_{V;\omega} = n+2$.
\end{example}

\begin{sloppypar}

\begin{example} Let $k=1$, $f_1=f$, $\omega=dx_1$. One has
\begin{eqnarray}
{\cal A}_{V;\omega} & = & {\cal O}_{\CC^n,0}/\left\langle f,
{\textstyle \frac{\partial f}{\partial x_2}, \ldots , \frac{\partial f}{\partial x_n}} \right\rangle,
\label{Equ1} \\
\widetilde{J} & = & {\rm Jac}\left(f, {\textstyle \frac{\partial f}{\partial x_2}, \ldots ,
\frac{\partial f}{\partial x_n}} \right) \left( {\textstyle \frac{\partial f}{\partial
x_1}} \right)^{2-n}.
\end{eqnarray}
One can consider the map
$$
\left(f, {\textstyle \frac{\partial f}{\partial x_2},\ldots,\frac{\partial f}{\partial x_n}} \right):
(\CC^n,0) \to (\CC^n,0).
$$
The linear function $R_{\ELKh}$ (and the corresponding (non-degenerate)
quadratic form $Q_{\ELKh}$) on the algebra ${\cal A}_{V;\omega}$ which
participate in the Eisenbud-Levine-Khimshiashvili theory can be defined
just by the formula (\ref{Formula1}) applied to the 1-form $fdx_1 +
\frac{\partial f}{\partial x_2}dx_2 + \ldots +
\frac{\partial f}{\partial x_n}dx_n$ on $(\CC^n,0)$ (with
$\widetilde{J}={\rm Jac}(f, \frac{\partial f}{\partial x_2}, \ldots ,
\frac{\partial f}{\partial x_n})$). With (\ref{Equ1}) this implies that
\begin{eqnarray*}
R_{V,\omega}(\varphi) & = & R_{\ELKh}(({\textstyle \frac{\partial f}{\partial x_1}})^{n-2}
\varphi), \\
Q^{\cal A}_{V;\omega}(\varphi_1,\varphi_2) & = & Q_{\ELKh}(({\textstyle\frac{\partial f}{\partial
x_1}})^{n-2} \varphi_1,\varphi_2).
\end{eqnarray*}
Therefore the rank of $Q^{\cal A}_{V;\omega}$ is equal to the rank of the
multiplication by $(\partial f/\partial
x_1)^{n-2}$ and
$$
Q^{\cal A}_{V;\omega} = 0 \Leftrightarrow  \left( {\textstyle \frac{\partial
f}{\partial x_1}} \right)^{n-2} \in \left\langle f, {\textstyle \frac{\partial f}{\partial
x_2}, \ldots , \frac{\partial f}{\partial x_n}} \right\rangle.
$$
In particular,
for $n=2$ (i.e.\ for plane curves) these two quadratic forms coincide.
\end{example}

\end{sloppypar}

For a germ of a reduced curve $C$ and for a meromorphic 1-form  $\omega$ on
it, Montaldi and van Straten \cite{MvS} introduced two quadratic forms
$\psi^\pm_\omega$ on certain modules $R^\pm(\omega)$. If the 1-form $\omega$ is
holomorphic, one has $R^-(\omega)=0$. If, in addition, the curve $C$ is an
{\icis} defined by a map $f:(\CC^n,0) \to (\CC^{n-1},0)$ then the module
$R^+(\omega)$ coincides with the algebra ${\cal A}_{C;\omega}$
(cf.\ \cite[(1.5)]{MvS}) and the quadratic form $\psi^+_\omega$ coincides
with $Q^{\cal A}_{V;\omega}$. (This can be derived,
e.g., from the fact that the
dimension of $R^+(\omega)$ satisfies the law of conservation of number
\cite[Theorem (1.7)]{MvS} and the quadratic forms $\psi^+_\omega$ and
$Q^{\cal A}_{V;\omega}$ coincide for a smooth curve.)


\section{A quadratic form on the module of differential forms}
There is a module which as a $\CC$-vector space has the same dimension as
${\cal A}_{V;\omega}$. This is the ${\cal O}_{V;\omega}$-module
$$\Omega_{V;\omega} := \Omega^{n-k}_{\CC^n,0}/(f_i \Omega^{n-k}_{\CC^n,0},
df_i \wedge \Omega^{n-k-1}_{\CC^n,0}, \omega \wedge
\Omega^{n-k-1}_{\CC^n,0}) (= \Omega^{n-k}_{V,0}/ \omega \wedge
\Omega^{n-k-1}_{V,0}).$$
The fact that $\dim_\CC \Omega_{V;\omega} = \dim_\CC {\cal A}_{V;\omega}$
is proved in \cite{Gr1}. (There it is proved for $\omega=df$, but
G.-M.~Greuel informed us that there is no difference for the general case.)

One can define a quadratic form $Q^\Omega_{V;\omega}$ on
$\Omega_{V;\omega}$ in the spirit of the definition above. Let $\eta_1,
\eta_2 \in \Omega_{V;\omega}$ and let $\eta_1$ and $\eta_2$ be
representatives of the elements $\eta_1$ and $\eta_2$ in
$\Omega^{n-k}_{\CC^n,0}$. We consider a deformation $F_\eps$ of the map
$f:(\CC^n,0) \to (\CC^k,0)$ with a deformation $\omega_\eps$ of the 1-form
$\omega$. For $\eps \not\in \Sigma$, on the smooth level set $\{
f_\eps=0\}$ the tensor product $\eta_1 \otimes \eta_2$ of the $(n-k)$-forms
$\eta_1$ and $\eta_2$ is a tensor of the same type as $
\widetilde{J}(dy^1 \wedge \ldots \wedge dy^{n-k})^2$. Since all the
singular points $P$ of the
1-form $\omega_\eps$ on the level set $\{f_\eps=0\}$ are non-degenerate, the latter
tensor (the Jacobian) does not tend to zero at these points. Now one can
define $Q^\Omega_{V;\omega}(\eta_1,\eta_2)$ as
\begin{equation}
Q^\Omega_{V;\omega}(\eta_1,\eta_2):= \lim_{\eps \to 0} \sum_{P \in {\rm
Sing}_{V_\eps} \omega_\eps} \frac{(\eta_1 \otimes \eta_2)_{|P}}{\widetilde{J}(dy^1
\wedge \ldots \wedge dy^{n-k})^{\otimes 2}}. \label{Formula3}
\end{equation}
It is possible to show (in the same way as above) that this definition
makes sense (i.e.\ that the limit exists) and that the result does not
depend on the choice of representatives of the elements $\eta_1, \eta_2
\in \Omega_{V;\omega}$. On the other hand, this follows directly from another
description of the quadratic form $Q^\Omega_{V;\omega}$.

Let $\Lambda : \Omega_{V;\omega} \to {\cal A}_{V;\omega}$ be the map which
sends an $(n-k)$-form $\eta \in \Omega^{n-k}_{\CC^n,0}$ to the function
$$\frac{df_1 \wedge \ldots \wedge df_k \wedge \eta}{dx_1 \wedge \ldots
\wedge dx_n}$$
(one can easily see that this mapping is well-defined). Now one has
$$Q^\Omega_{V;\omega}(\eta_1,\eta_2)=Q^{\cal A}_{V;\omega}(\Lambda \eta_1,
\Lambda \eta_2).$$

\begin{remarks} 1. Again, if the map $f$ and the 1-form $\omega$ are real,
the quadratic form $Q^\Omega_{V;\omega}$ is real as well.

2. We have indicated that the linear function $R$ and therefore the
quadratic form  $Q^{\cal A}_{V;\omega}$ on ${\cal A}_{V;\omega}$ are
defined only up to multiplication by a unit in the algebra ${\cal
A}_{V;\omega}$. The mapping $\Lambda$ is also defined up to multiplication
by a unit (in fact the same one). However, the quadratic form
$Q^\Omega_{V;\omega}$ is absolutely well-defined. This follows from the
representation (\ref{Formula3}) where the denominator is a well-defined
tensor at the singular points $P$.
\end{remarks}

\begin{corollary} ${\rm rk}\, Q^\Omega_{V;\omega} \leq {\rm rk}\, Q^{\cal
A}_{V;\omega}$.
\end{corollary}

Let us show that the difference ${\rm rk}\, Q^{\cal A}_{V;\omega} - {\rm
rk}\, Q^\Omega_{V;\omega}$   is bounded from above by a constant which does
not depend on the 1-form $\omega$.

Let $T\Omega^{n-k}_{V,0}$ be the torsion module of the module
$$\Omega^{n-k}_{V,0} = \Omega^{n-k}_{\CC^n,0}/(f_i \Omega^{n-k}_{\CC^n,0},
df_i \wedge \Omega^{n-k-1}_{\CC^n,0})$$
of differentiable $(n-k)$-forms on the {\icis} $(V,0)$. The dimension of
$T\Omega^{n-k}_{V,0}$ as a $\CC$-vector space is denoted by $\tau'$ (cf.\
\cite{Gr2}). (This dimension is finite since the corresponding sheaf is a
coherent one concentrated at the origin.)

\begin{proposition}  \label{Prop3}
One has
$${\rm Ker}\, \Lambda \cong T\Omega^{n-k}_{V,0}$$
and therefore $\dim_\CC {\rm Ker}\, \Lambda = \tau'$.
\end{proposition}

\begin{proof} There exists the following commutative diagram with exact
rows:
$$\diagram
0 \rto & \omega \wedge \Omega^{n-k-1}_{V,0} \rto
\dto_{\lambda'} & \Omega^{n-k}_{V,0} \rto \dto_{\lambda=\frac{df \wedge}{dx}}
& \Omega_{V;\omega} \rto \dto_{\Lambda} & 0 \\
0 \rto & \langle m_I \rangle  \rto &
{\cal O}_{V,0}  \rto &
{\cal A}_{V;\omega} \rto & 0
\enddiagram
$$
Here $m_I$ is the $(k+1) \times (k+1)$-minor of the matrix
$$\begin{pmatrix} \frac{\partial f_1}{\partial x_1} & {\cdots } &
 \frac{\partial f_1}{\partial x_{n}} \\ {\vdots} & {\cdots} &
{\vdots} \\ \frac{\partial
f_k}{\partial x_1} & \cdots &
 \frac{\partial f_k}{\partial x_{n}} \\
A_1 & \cdots & A_{n}
\end{pmatrix}$$
consisting of columns with indices from $I \subset \{1, \ldots ,n\}$,
$\langle m_I \rangle \subset {\cal O}_{V,0}$ is the ideal generated by the minors $m_I$. The
restriction $\lambda' : \omega \wedge \Omega^{n-k-1}_{V,0} \to \langle m_I \rangle$
of the mapping $\lambda= \frac{df \wedge}{dx}$ to the submodule
$\omega \wedge \Omega^{n-k-1}_{V,0}$ is surjective (since $\omega\wedge\wedge_{j \not\in I} dx_j$ 
maps to $\pm m_I$). One has
$\dim_\CC
\Omega_{V;\omega} = \dim_\CC {\cal A}_{V;\omega}$ and therefore
$\dim_\CC {\rm Ker}\, \Lambda  =  \dim_\CC {\rm Coker}\, \Lambda$.
Moreover, one has
$$
{\rm Coker}\, \lambda  =  {\cal O}_{V,0}/\langle \det \left( {\textstyle \frac{\partial f_i}{\partial
x_{s_j}}} \right), 0 \leq s_1 < \ldots s_k \leq n \rangle.
$$
This implies that $\dim_\CC {\rm Coker}\, \lambda  =  \tau'$  
(cf.\ \cite[0.1 Satz]{Gr2}). On the other hand
${\rm Ker}\, \lambda = T\Omega^{n-k}_{V,0}$ \cite[Proof of Proposition~1.11]{Gr1}.
Now the Snake lemma (see, e.g., \cite[Exercise
A3.10]{Eis}) yields the statement.
\end{proof}

\begin{corollary}
One has
$${\rm cork}\, Q^\Omega_{V;\omega} \geq \tau', \quad 0 \leq {\rm rk}\,
Q^{\cal A}_{V;\omega} - {\rm rk}\, Q^\Omega_{V;\omega} \leq 2\tau'.$$
\end{corollary}

\begin{examples}

1. In the situation of Example~1, the module $\Omega_{V;\omega}$ (as a $\CC$-vector space) is generated
by the forms
\begin{eqnarray*}
\alpha_i & := & dx_1 \wedge \ldots \wedge \widehat{dx_i} \wedge \ldots \wedge dx_n,\\
\beta_i & := & x_i dx_1 \wedge \ldots \wedge \widehat{dx_i} \wedge \ldots \wedge dx_n
\end{eqnarray*}
($i=1, \ldots , n$). We have $\Lambda(\alpha_i)=\pm x_i$, $\Lambda(\beta_i)=\pm x_i^2$ (recall that
$\sum x_i^2 =0$ in ${\cal A}_{V;\omega}$). Therefore
$$
Q^\Omega_{V;\omega}(\alpha_i,\alpha_i) = \frac{2}{\prod_{j \neq i} (a_j - a_i)}
$$
and $Q^\Omega_{V;\omega}(\eta_1,\eta_2)=0$ for all other pairs of the generators listed above. Hence
${\rm rk}\, Q^\Omega_{V;\omega} =n$.

2. In Example~2, there is a natural isomorphism between
$\Omega_{V;\omega}$ and ${\cal A}_{V;\omega}$ (as ${\cal O}_{\CC^n,0}$-modules). Under this
isomorphism the map $\Lambda$ becomes the multiplication by $\partial
f/\partial
x_1$. Therefore
$$
Q^\Omega_{V;\omega}(\varphi_1,\varphi_2) = Q_{\ELKh}(({\textstyle \frac{\partial f}{\partial
x_1}})^{n-1} \varphi_1,\varphi_2),
$$
the rank of $Q^\Omega_{V;\omega}$ is equal to the rank of the
multiplication by $(\partial f/\partial
x_1)^{n-1}$ and
$$
Q^\Omega_{V;\omega} = 0 \Leftrightarrow  \left( {\textstyle \frac{\partial f}{\partial
x_1}} \right)^{n-1} \in \left\langle f,
{\textstyle \frac{\partial f}{\partial x_2}, \ldots , \frac{\partial f}{\partial x_n}}
\right\rangle.
$$
\end{examples}

\bigskip
\noindent Universit\"{a}t Hannover, Institut f\"{u}r Mathematik \\
Postfach 6009, D-30060 Hannover, Germany \\
E-mail: ebeling@math.uni-hannover.de\\

\medskip
\noindent Moscow State University, Faculty of Mechanics and Mathematics\\
Moscow, 119992, Russia\\
E-mail: sabir@mccme.ru


\begin{thebibliography}{GM2}

\bibitem[AGV]{AGV} Arnold, V.~I., Gusein-Zade, S.~M., Varchenko, A.~N.:
Singularities of Differentiable Maps, Vol. I, Birkh\"auser,
Boston Basel Berlin, 1985.

\bibitem[EG1]{EG1} Ebeling, W., Gusein-Zade, S.~M.: Indices of 1-forms
on an isolated complete intersection singularity. Mosc.
Math. J. {\bf 3} (2003), no.2, 439--455.

\bibitem[EG2]{EG2} Ebeling, W., Gusein-Zade, S.~M.: Indices of vector fields or
1-forms and characteristic numbers. math.AG/0303330, Bull. London Math. Soc.
(to appear).

\bibitem[Ei]{Eis} Eisenbud, D.: Commutative Algebra with a View Toward
Algebraic Geometry. Graduate texts in Math. {\bf 150}, Springer Verlag, 
1994.

\bibitem[EL]{EL} Eisenbud D., Levine, H.: An algebraic
formula for the degree of a $C^\infty$ map germ. Ann. Math. {\bf 106},
19--38 (1977).

\bibitem[GM1]{GMM1} G\'omez-Mont, X., Marde\v{s}i\'c, P.: The index of a
vector field
tangent to a hypersurface and the signature of the relative Jacobian
determinant.
Ann. Inst. Fourier {\bf 47}, 1523--1539 (1997).

\bibitem[GM2]{GMM2} G\'omez-Mont, X., Marde\v{s}i\'c, P.:
The index of a vector field tangent to an odd-dimensional
hypersurface and the signature of the relative Hessian.
Funktsional. Anal. i Prilozhen. {\bf 33}, no.1, 1--13 (1999)
(English translation in Funct. Anal. and Applications,
{\bf 33}, no.1, 1--10 (1999)).

\bibitem[Gr1]{Gr1}
Greuel, G.-M.: Der Gau\ss-Manin-Zusammenhang
isolierter Singularit\"aten von vollst\"andigen Durchschnitten.
Math. Ann. {\bf 214}, 235--266 (1975).

\bibitem[Gr2]{Gr2}
Greuel, G.-M.: Dualit\"at in der lokalen Kohomologie isolierter
Singularit\"aten. Math. Ann. {\bf 250}, 157--173 (1980).

\bibitem[Kh]{Kh} Khimshiashvili, G.~N.: On the local degree
of a smooth map. Comm. Acad. Sci. Georgian SSR. {\bf 85}, no.2,
309--311 (1977) (in Russian).

\bibitem[Le]{Le} L\^e, D.~T.: Computation of the Milnor number
of an isolated singularity of a complete intersection.
Funktsional. Anal. i Prilozhen. {\bf 8}, no.2,  45--49
(1974) (English translation in Funct. Anal. and Applications, {\bf 8}, no.2, 127--131 (1974)).

\bibitem[MvS]{MvS} Montaldi, J., van Straten, D.: One-forms on singular curves
and the topology of real curve singularities. Topology {\bf 29}, 501--510
(1990).

\bibitem[Sz]{Sz} Szafraniec, Z.: On topological invariants of
real analytic
singularities. Math. Proc. Camb. Phil. Soc. {\bf 130}, 13--24 (2001).

\end{thebibliography}
\end{document}